%&amstex

\documentstyle{amsppt}

\loadbold
\mag=1200
\refstyle{A}
\pagewidth{32pc}

% ---------------------------------------
% -------- Definitions ------------

% ------ Fonts ----------
\font\tiny=cmmi6
\font\verysmall=cmr6

% ------ Symbols -------
\def\f{\goth F}
\def\g{\goth G}
\def\fr{\goth F_{\text{\tiny R}}}
\def\m{\goth m}
\def\p{\goth p}
\def\I{\goth I}
\def\tr{\text{\rm Tr}}
\def\pr{\text{\rm pr}}
\def\can{\text{\rm can}}
\def\id{\text{\rm id}}
\def\spec{\text{\rm Spec}}
\def\mor{\text{\rm Mor}}

\def\A{\Cal A}
\def\B{\Cal B}
\def\O{\Cal O}
\def\Ot{\tilde{\O}}
\def\R{\Cal R}

\def\nrd{\text{\rm Nrd}}
\def\chr{\text{\rm chr}}
\def\n{\text{\rm N}}
\def\SL{\text{\rm SL}}
\def\GL{\text{\rm GL}}
\def\SU{\text{\rm SU}}
\def\U{\text{\rm U}}

\def\Or{\text{\rm O}}
\def\SO{\text{\rm SO}}

\def\End{\text{\rm End}}
\def\N{\Bbb N}
\def\s{\sigma}
\def\st{\tilde{\s}}

\def\hml{H_{\text{\verysmall  \'{e}t}}^1}
\def\ainv{{\A^*}^{1-\s}}

\def\e{\varepsilon}
\def\et{\tilde{\e}}
\def\St{\tilde{S}}
\def\a{\alpha}

%-------- Arrows ---------
\def\ra{\rightarrow}
\def\ar{\longrightarrow}
\def\xra#1{\overset{#1}\to{\rightarrow}}
\def\hra#1{\overset{#1}\to{\hookrightarrow}}
\def\lra#1{\!\!\!\! \overset{#1}\to{\ar} \!\!\!\!}
\def\da{\downarrow}

% -------- Stars -----------
\def\ss{\hbox{($**$)}}
\def\ssp{\hbox{($**'$)}}
\def\sss{\hbox{($*\!*\!*$)}}

\def\Ap{\A_{_{+}}}
\def\Am{\A_{_{-}}}

%------------------------------------------
%------------------------------------------
\topmatter
\title 
On Grothendieck's Conjecture 
about Principal Homogeneous Spaces\\
for Some Classical Algebraic Groups
\endtitle

\author
Kirill Zainoulline
\endauthor

\affil
Laboratory of Algebra, \\
Steklov Mathematical Institute, \\
St.~Petersburg Department, \\
Russian Academy of Sciences, \\
Fontanka~27, St.~Petersburg 191011, Russia \\
E-mail: kirill\@pdmi.ras.ru \\
URL: http://www.pdmi.ras.ru/\~{}kirill/
\endaffil

\date 
February 4, 1999
\enddate

\abstract
In the present paper we investigate the question about
the injectivity of the map $\f(R) \ra \f(K)$ induced by the
canonical inclusion of a local regular ring of geometric
type $R$ to its field of fractions $K$  for a homotopy
invariant functor $\f$ with transfers satisfying some additional properties. 
As an application we get the original proof of Special Unitary Case
of Grothendieck's conjecture about principal homogeneous spaces and 
some other interesting examples.
\endabstract
\endtopmatter

\rightheadtext{Kirill Zainoulline}
\leftheadtext{On Grothendieck's Conjecture}

\document

Let $R$ be a local regular ring of geometric type, 
i.e. $R$ is the local ring of a point
of a smooth affine variety over a field $k$.
Let $K$ be its field of fractions. 
Let $\f$ be a covariant functor from the category of $R$-algebras 
to the category of abelian groups.

Assume that $\f$ is a homotopy invariant functor with transfers satisfying
some additional properties.
Our main goal in the present article is to prove the injectivity of the map
$\f(R)\ra \f(K)$ induced by the canonical inclusion (sections~1 and 2).

As a consequence of this result we get the positive solution 
of Special Linear Case of Grothendieck's conjecture 
about principal homogeneous spaces \cite{Gr} and
some other important applications (section~3).
In particular, using the main result of \cite{PO2} and 
Norm Principle for the unitary group (section 4)
we get the original proof of Special Unitary Case.
 
Observe that Special Linear Case was originally proved in the work 
of I.~Panin and A.~Suslin \cite{PS} but the method they have used 
doesn't work for the other interesting cases of Grothendieck's conjecture,
say for Special Unitary Case.

Finally, using a well-known theorem of D.~Popescu \cite{PO1}
we generalize some our results to the case when
$R$ is a local regular ring containing a field (section 5).

Our work was motivated by the paper of V.~Voevodsky \cite{Vo} and, 
mainly, by the paper of I.~Panin and M.~Ojanguren
on one Grothendieck's conjecture for hermitian spaces \cite{PO2}.
The point was to offer a good axiomatization for the method
they introduced.   

\subhead{Acknowledgements} \endsubhead
I~am personally grateful to Ivan Panin for his attention and help. 
I~thank University Franche-Comte Besancon and Eva Bayer-Fluckiger 
for hospitality and useful discussions on the subject of this article.
I~am very grateful to Vladimir Chernousov 
for the close reading of my paper during preparation and
for the ideas led to the proof of Norm Principle for the unitary group.

\head{Agreements}\endhead

All rings are assumed to be commutative noetherian with unit.
By $k$ we will always mean an infinite field. 

Let $A$ and $S$ be any rings.
By an $A$-algebra $S$ we will mean the pair $(S, i)$, 
where $i: A \ra S$ is a ring homomorphism. 
Sometimes, we will write just $S$, instead of the pair $(S,i)$, 
keeping in mind the homomorphism $i$.
We will denote an $A$-algebra $S$ by $A \xra{i} S$.

By a morphism $f : (S, i) \ra (S',i')$ between two $A$-algebras 
we will mean the commutative diagram:
$$
\CD
A @>{i}>> S \\
@|    @VV{f}V \\
A @>>{i'}> S'
\endCD
$$

Let $\f$ be a covariant functor on the category of $A$-algebras 
to the category of abelian groups.
Let $A \xra{i} R$ be an $A$-algebra.
By restriction of $\f$ to the category of $R$-algebras along $i$ 
we will call the functor denoted by $\fr$ and 
given as follows:
$\fr(R \xra{t} T)=
\f(A \xra{i} R \xra{t} T)$
on objects and 
$$
\mor(\fr(R \xra{t_1} T_1),\fr(R \xra{t_2} T_2))=
\mor(\f(A \xra{i} R \xra{t_1} T_1),\f(A \xra{i} R \xra{t_2} T_2 )),
$$
on morphisms for any $R$-algebra $T$, $T_1$, $T_2$.

Further in the paper we will use the result of Grothendieck
(\cite{Ei}, Corollary 18.17)  which says that if we have a finite extension 
$A \hra{} B$ of essentially smooth $k$-algebras  then $B$ is  finitely
generated projective as the $A$-module. 

% ***************************************************************
% ***************************************************************
% ***************************************************************
\head{1. Constant Case}\endhead

Let $A$ be a smooth $k$-algebra.
Let $R=A_\p$ be the local ring at a prime ideal $\p$. 
By $\m_R$ we will denote the corresponding maximal ideal of $R$.
Thus, we have the $A$-algebra $A \hra{i_R} R$, 
where $i_R$ is the canonical inclusion.

Let $\f$ be a covariant functor from the category of $k$-algebras 
to the category of abelian groups.
By $\fr$ we denote its restriction to the category of $R$-algebras 
along $k \hra{} A \hra{i_R} R$. 
Let functors $\f$ and $\fr$ satisfy the following list of axioms:

\subsubhead{Axiom for the functor $\f$}\endsubsubhead

\itemitem{\bf C.}(continuity)
For any $A$-algebra $S$ essentially smooth over $k$ and
for any multiplicative system $M$ in $S$ the canonical map
$\varinjlim_{g \in M} \f(S_g) \ra \f(M^{-1}S)$ is an isomorphism,
where $M^{-1}S$ is the localization of $S$ with respect to $M$.

\subsubhead{Axioms for the functor $\fr$}\endsubsubhead

For any $R$-algebra $T$ finitely generated and projective as the $R$-module

\itemitem{\bf TE.} (existance)
It is given a homomorphism $\tr_R^T: \fr(T) \ra \fr(R)$ called transfer map;

\itemitem{\bf TA.} (additivity)
For every element $x \in \fr(R\times T)$ with $x_R=\pr_R^*(x)\in \fr(R)$ and
$x_T=\pr_T^*(x)\in \fr(T)$ the relation
$\tr_R^{R\times T}(x)= x_R + \tr_R^T(x_T)$ holds in $\fr(R)$,
where $\pr_R^*$ and $\pr_T^*$ are induced by projections;

\itemitem{\bf TB.} (base changing and homotopy invariance)
For any $R[t]$-algebra $S$ finitely generated projective as the $R[t]$-module
(thus, $S/(t)$ and $S/(t-1)$ are finitely generated projective as the $R$-modules)
the following diagram commutes:
$$
\CD
\fr(S) @>{\can_0^*}>> \fr(S/(t)) \\
@V{\can_1^*}VV        @VV{\tr_0}V \\
\fr(S/(t-1)) @>>{\tr_1}> \fr(R)
\endCD
$$
where $\can_0^*$, $\can_1^*$ are induced by the canonical projections and
$\tr_0$, $\tr_1$ denote the corresponding transfer maps
$\tr_R^{S/(t)}$ and $\tr_R^{S/(t-1)}$.

Let $K$ be a field of fractions of the ring $A$. 
Then our aim is to prove: 

\proclaim{Theorem  {\rm (Constant case)}}
Let $\f$ be the functor on the category of $k$-algebras and 
let $\fr$ be its restriction to the category of $R$-algebras.
If  $\f$ and $\fr$ satisfy axioms {\bf C}, {\bf TE}, {\bf TA}, {\bf TB} 
then the homomorphism $\f(R) \ra \f(K)$ 
induced by the canonical inclusion is injective.
\endproclaim

\subhead {Remark~1} \endsubhead
To get a better feeling for the axioms above observe that 
{\bf TE}, {\bf TA}, {\bf TB}
are the consequences of the following more strong conditions:

For any  $R$-algebra $S$ essentially smooth over $k$ and 
for any $S$-algebras $T_1$, $T_2$ and $T$ finitely generated projective as the $S$-modules

\itemitem{\bf H.} (homotopy invariance)
The map $\fr(S) \ra \fr(S[t])$ induced by the inclusion is an isomorphism;

\itemitem{\bf TE$'$.} (existance) 
It is given a homomorphism $\tr_S^T: \fr(T) \ra \fr(S)$ called transfer map; 

\itemitem{\bf TA$'$.} (additivity) 
Let $T=T_1\times T_2$. 
For every $x \in \fr(T)$, $x_1=\pr_1^*(x)\in \fr(T_1)$ and 
$x_2=\pr_2^*(x)\in \fr(T_2)$ the relation
$\tr_S^T(x)= \tr_S^{T_1}(x_1) + \tr_S^{T_2}(x_2)$ holds in $\fr(S)$, 
where $\pr_i^*: \fr(T_1\times T_2) \ra \fr(T_i)$ are
induced by projections; 

\itemitem{\bf TB$'$.} (base changing) 
For any extension $U/S$ of an $R$-algebras with $U$ essentially smooth over $k$ the following
diagram commutes: 
$$
\CD
\fr(T) @>>> \fr(U\otimes_S T) \\
@V{\tr_S^T}VV  @VV{\tr_U^{U\otimes_ST}}V \\
\fr(S) @>>> \fr(U)
\endCD
$$

Indeed, to show that the axiom {\bf TB} 
is a consequence of the axioms {\bf H} and {\bf TB$'$}
let look at two commutative diagrams arising 
when we  apply {\bf TB$'$} to the extension $S/R[t]$ and 
evaluations $i_0: R[t] \xra{t\mapsto 0} R$, 
$i_1: R[t]  \xra{t\mapsto 1} R$ at $0$
and $1$ correspondingly: 
$$
\CD
\fr(S) @>{\can_0^*}>> \fr(S/(t)) \\
@V{\tr^S_{R[t]}}VV        @VV{\tr_0}V \\
\fr(R[t]) @>>{i_0^*}> \fr(R)
\endCD 
\qquad \qquad
\CD
\fr(S) @>{\can_1^*}>> \fr(S/(t-1)) \\
@V{\tr^S_{R[t]}}VV        @VV{\tr_1}V \\
\fr(R[t]) @>>{i_1^*}> \fr(R)
\endCD 
$$

By axiom {\bf H}  the map $\fr(R)
\ra \fr(R[t])$  
induced by the inclusion is an isomorphism. 
Since the compositions 
$R \hra{} R[t] \xra{i_0} R$ and $R \hra{} R[t] \xra{i_1} R$
are identities, evaluations $i_0^*$ and $i_1^*$ coinside.

Gluing together these two diagrams via the left-down part
we get the required one from {\bf TB}. \qed

The rest of this section is organized as follows: 
1.1 contains the proof of Specialization Lemma,
1.2 is devoted to a version of Quillen's Trick,
and the last subsection~1.3 contains the proof of our theorem.

\subhead\nofrills{1.1. Specialization Lemma}\endsubhead

Let $R$ be a local regular ring of a smooth $k$-algebra $A$ and 
let $R \xra{i} S$ be a given $R$-algebra with an element $f \in S$
such that:

\itemitem{\bf S1.}
a)~$S$ is finite over $R[t]$ for some specially chosen $t \in S$;\ \ 
b)~The quotient $S/(f)$ is finite over $R$;

\itemitem{\bf S2.}
There is an augmentation map $\e: S \ra R$
such that the composition $R \xra{i} S \xra{\e} R$ is the identity;

\itemitem{\bf S3.}
a)~$S$ is essentially smooth over the field $k$;\ \ 
b)~$S/\m_RS$ is smooth over  the residue field $R/\m_R$
at the maximal ideal $\e^{-1}(\m_R)/\m_RS$.

\subhead {\bf Remark~2} \endsubhead
To see that {\bf S3.(b)} is described correctly observe that we have
an obvious inclusion $\m_RS\subset \e^{-1}(\m_R)$.
Since $R$ is a local we conclude that the ideal $\e^{-1}(\m_R)$ 
is the unique maximal ideal lying over the kernel of
augmentation $I=\ker(\e)$. \qed

Now we are going to prove:
\proclaim{Theorem {\rm (Specialization Lemma)}}
Let $(R,S,f)$ be the triple satisfying conditions {\bf S1}, {\bf S2} and  {\bf S3}.
Let $\g$ be a covariant functor from the category of $R$-algebras 
to the category of abelian groups 
satisfying axioms {\bf TE}, {\bf TA} and {\bf TB}.

If the image in $\g(S_f)$ of an element $\a_S \in \g(S)$ is trivial 
then $\e^*(\a_S)=0$ in $\g(R)$.
\endproclaim

\demo{Proof}
To simplify the notation, for any element $\a_S \in \g(S)$,
we will denote by $\a_{S_f}$ its image under the map
induced by the inclusion $S \hookrightarrow S_f$ and
by $\a_0$, $\a_1$ and $\a_I$  --
its images under the maps induced by the canonical projections
$\can_0^*$, $\can_1^*$ and $\can_I^*$ for the ideal $I$ correspondingly.

Since {\bf S1} and {\bf S3} hold for the triple $(R,S,f)$, by Geometric Presentation Lemma
(\cite{PO2}, Lemma 10.1)
for the given $f \in S$  and $t \in S$ we can find an element  $t'\in S$ such that:

\itemitem{\bf G1.} $S$ is finite over $R[t']$;

\itemitem{\bf G2.} There exists an ideal $J$ coprime with $I$ and with the property $I\cap J=(t')$; 

\itemitem{\bf G3.} $(f)$ is coprime with $J$ and with $(t'-1)$.

\enddemo

\subhead{I}\endsubhead
The last property {\bf G3} means that we can factorize the
canonical projections $S \xra{\can_J} S/J$ and 
$S \xra{\can_1} S/(t'-1)$ through the localization $S_f$.

Hence, if $\a_{S_f}=0$ in $\g(S_f)$ for some $\a_S \in \g(S)$  
then $\a_J=\can_J^*(\a_S)=0$ and  $\a_1=\can_1^*(\a_S)=0$.

\subhead{II} \endsubhead
By the theorem of Grothendieck
(see the end of Agreements), 
since $S$ is essentially smooth over $k$ and $S$ is finite over $R[t']$ by {\bf G1}, 
$S$ is finitely generated projective over $R[t']$.

Now look at the commutative diagram from {\bf TB}:
$$
\CD
\g(S) @>{\can_0^*}>> \g(S/(t')) \\
@V{\can_1^*}VV        @VV{\tr_0}V \\
\g(S/(t'-1)) @>>{\tr_1}> \g(R)
\endCD
$$
By the previous step we get $\tr_0(\a_0)=\tr_1(\a_1)=0$.

\subhead{III}\endsubhead
By {\bf G2}, there is the decomposition
$S/(t') \cong R \times S/J$ via the projections
$S/(t') \xra{\can_{0I}} S/I
\xra{\bar{\e}} R$ and
$S/(t') \xra{\can_{0J}} S/J$, where $\bar{\e}$ is an isomorphism 
taken from the obvious decomposition $\e: S \xra{\can_I} S/I \xra{\bar{\e}} R$.

Since $S/(t')$ is finitely generated projective over $R$,
the quotient  $S/J$ is finitely generated
projective over $R$, too. Therefore, by {\bf TE} we have well-defined
transfer map $\tr_J: \g(S/J) \ra \g(R)$.

By the additivity {\bf TA} of the transfer  and by the second step, we can write:
$$
0=\tr_0(\a_0)=
\bar{\e}^*(\can_{0I}^*(\a_0))+ \tr_J(\can_{0J}^*(\a_0)).
$$

Since $\can_{0I}^*(\a_0)=\a_I$ and $\can_{0J}^*(\a_0)=\a_J$,
we can rewrite the last relation as:
$$
0=\bar{\e}^*(\a_I)+ \tr_J(\a_J).
$$

By the first step $\a_J=0$, thus, one gets $\e^*(\a_S)=\bar{\e}^*(\a_I)=0$.
And we have done. \qed

\subhead\nofrills{1.2. A version of Quillen's Trick}\endsubhead

Let $(A, R, f)$ be a triple, where
$A$ is a smooth $d$-dimensional $k$-algebra,
$R$ be a local regular ring at the prime ideal $\p$ of $A$ and
$f \in \p$ be some fixed regular element.

We would like to produce certain extension $S$ of the ring $R$, 
such that properties {\bf S1}, {\bf S2}, {\bf S3} 
of subsection~1.1 are satisfied.

We need the following lemma of Quillen (see \cite{Qu}, Lemma 5.12):

\proclaim{Theorem {\rm (Quillen's Lemma)}}
Let $A$ be a smooth finite type algebra of dimension $d$ over a field $k$,
let $f$ be a regular element of $A$,
and let $\I$ be a finite subset of $\spec A$.
Then there exist elements $x_1,\ldots,x_d$ in $A$
algebraically independed over $k$
and such that  if $P=k[x_1,\ldots,x_{d-1}] \hra{q} A$, then \ \ 
{\rm i)} A/(f) is finite over $P$; \ \ 
{\rm ii)} A is smooth over $P$ at points of $\I$; \  and \ 
{\rm iii)} the inclusion $q$ factors as
$q: P \hra{} P[x_d] \xra{q_1} A$, where $q_1$ is finite.
\endproclaim

Set $\I=\{\p\}$ and apply Quillen's Lemma to the given
pair $(A,f)$. 

Look at the canonical tensor product diagram 
(we keep all notation coming from Quillen's Lemma):
$$
\CD
A @>{i_S}>> A\otimes_P R \\
@A{q}AA  @AA{i}A \\
P @>>{r}> R
\endCD
$$
where the map $r$ is the composition 
$r: P \hra{q} A \hra{i_R} R$ and 
$i_S: a \mapsto a\otimes 1$, $i: r \mapsto 1\otimes r$. 

We are going to show that the triple 
$(R,S, f\otimes 1)$ with $S=A\otimes_P R$ satisfies properties {\bf S1}, {\bf S2}, {\bf S3}
of subsection~1.1:

\itemitem{\bf S1.} a) Since $A$ is finite over $P[x_d]$ via $q_1$, 
we get that  $A\otimes_P R$ is finite over $P[x_d]\otimes_P R=R[x_d]$; \ \ 
b) Obviously, $S/(f\otimes 1)=A/(f)\otimes_P R$ is finite over $R$.
Further we will identify $f$ with $f \otimes 1$ via the map $i_S$.

\itemitem{\bf S2.} Take the multiplication map $a\otimes r \mapsto ar$ as augmentation $\e$.

\itemitem{\bf S3.}  a) The fact that $A$ is smooth over the ring $P$ at $\p$ 
is equivalent to the smoothness of the map 
$r:P \hra{q} A \hra{i_R} A_\p$.
Hence,  $r$ is smooth. Thus, $i_S$ is smooth also and 
since $A$ is smooth over $k$, so is $S$; \ \ 
b) By the property of smooth extensions, 
since $A$ is smooth over the $P$ at $\p$,
we have that $S$ must be smooth over $R$ at all primes lying above $\p S$.
Since $\e^{-1}(\m_R)$ lies over $\p S$, $S$ is smooth over $R$ at $\e^{-1}(\m_R)$.

\subhead\nofrills{1.3. The Proof of Constant Case}\endsubhead

Let $\f$ be the functor on the category  of $k$-algebras and
let $\fr$ be its restriction to the category of $R$-algebras.
By the hypothesis of the Theorem $\f$ and $\fr$ satisfy  axioms
{\bf C},  {\bf TE}, {\bf TA} and {\bf TB}.
We have to show that the map $\f(R) \ra \f(K)$
induced by the canonical inclusion is injective.

Let $\a' \in \f(R)$ be such that its image $\a_K'$ in
$\f(K)$ is trivial.

Since $\f$ is continuous ({\bf C}),
we may assume that $\a'$ came from an element $\a$ in $\f(A_g)$,
where $A_g$ is the localization of $A$ at some $g \in A \setminus \p$,
and the image of $\a$ in $\f(K)$ is trivial.
Thus, we can write $\a'=i_R^*(\a)$, where $i_R: A_g \hra{} R$ is
the canonical inclusion.

Observe that $A_g$ is again smooth $k$-algebra and $R$ is its local regular ring,
therefore, we can replace $A$ by $A_g$ and
consider the element $\a$ as lying in $\f(A)$.

Using {\bf C} again we get that there exists an element
$f \in \p$ such that the image $\a_f$ in $\f(A_f)$ of the element $\a$ is trivial.

Hence, our problem has reduced to the following one:
\proclaim{Proposition}
For a given  $\a \in \f(A)$ and $f \in \p$
such that the image $\a_f$ is trivial in $\f(A_f)$
the element $\a'=i_R^*(\a)$ is trivial in $\f(R)$.
\endproclaim

\demo{Proof} 
The proof consists of two steps. On the first step for any triple $(A,R,f)$  
we build up a $3\times 3$ commutative diagram. 
On the second step by playing with this diagram for the specially chosen triple
we will complete the proof.
\enddemo

\subhead{I}\endsubhead
Let $R \xra{i} S$ be the extension of the ring $R$ built up by using Quillen's Trick for any triple
$(A,R,f)$ with $R=A_\p$ and regular $f \in \p$. 
In particular, there is the  commutative diagram:
$$
\CD
A @>{i_S}>> S \\
@| @VV{\e}V \\
A @>>{i_R}> R
\endCD
\tag{$*$}
$$
where $\e : S \ra R$ is the augmentation from {\bf S2} for the inclusion $i: R \ra S$.

Let $\fr$ be the restriction of the functor $\f$ to the
category of $R$-algebras.

By the commutativity of the diagram
$$
\CD
k @>>> A @>{i_S}>> S \\
@| @. @AA{i}A \\
k @>>> A @>>{i_R}> R
\endCD
$$
we have the identities
$$
\f(S)=\f(k \hra{} A \xra{i_S} S)=\f(k\hra{} A \hra{i_R} R  \xra{i} S)=\fr(R \xra{i} S)=\fr(S),
$$
$\f(S_f)=\fr(S_f)$ and $\f(R)=\fr(R)$.
Moreover, there is a natural identification of the functors $\f$
and $\fr$ restricted to the category
of $S$-algebras and the following diagram commutes:
$$
\CD
\f(S) @= \fr(S) \\
@V{\e^*}VV @VV{\e^*}V \\
\f(R) @= \fr(R)
\endCD
\tag{$**$}
$$

Consider the $3 \times 3$ commutative diagram
$$
\CD
\a_f \in @. \f(A_f) @>>> \f(S_f) @= \fr(S_f) @. \ni \a_{S_f}\\
@. @AAA @AAA @AAA @.\\
\a \in @. \f(A) @>{i_S^*}>> \f(S) @= \fr(S) @. \ni \a_S\\
@. @| @V{\e^*}VV @VV{\e^*}V @. \\
@. \f(A) @>>{i_R^*}> \f(R) @= \fr(R) @. \ni \a'
\endCD
\tag{$*\!*\!*$}
$$
where 
the upper squares are induced by the localization,
the left-down square is induced by ($*$) and
the right-down square coinsides with \ss.

\subhead{II}\endsubhead
Now let $\a \in \f(A)$ and $f$ be an element from $\p$ such
that $\a_f=0$.

Consider the commutative diagram  \sss{}  for the given triple $(A,R,f)$
and consider the elements 
$\a' \in \fr(R)$, $\a_S \in \fr(S)$, $\a_f \in \f(A_f)$, $\a_{S_f} \in \fr(S_f)$, 
where $\a'=i_R^*(\a)$, $\a_S=i_S^*(\a)$ and
$\a_f$, $\a_{S_f}$ are the images of the elements $\a$, $\a_S$ under the maps induced
by the canonical inclusions.

To finish our proof we apply Specialization
Lemma (see~1.1) to the right column of our diagram.

By the very construction the triple $(R,S,f)$ satisfies properties 
{\bf S1}, {\bf S2} and {\bf S3} of subsection~1.1.
And by the very assumption the functor
$\fr$ satisfies  axioms  {\bf TE}, {\bf TA}, {\bf TB}.
So we are under the hypothesis of Specialization Lemma.

Since $\a_f=0$ in $\f(A_f)$, we get $\a_{S_f}=0$ in $\fr(S_f)$.
Using Specialization Lemma we conclude that $\e^*(\a_S)=0$.
By the commutativity of the diagram $\e^*(\a_S)$ coincides with $\a'$, thus, $\a'=0$.
And we have finished. \qed

% ***************************************************************
% ***************************************************************
% ***************************************************************
\head{2. Non-Constant Case}\endhead

Consider more general situation, namely, 
let the functors $\f$ and $\fr$ (see section~1) are defined only on the category
of $A$-algebras and satisfy axioms {\bf C}, {\bf TE}, {\bf TA}, {\bf TB}.
And we still want to prove  the injectivity of the
homomorphism $\f(R) \ra \f(K)$. Recall that $A$ is smooth over $k$, $R=A_\p$ is the local
regular ring and $K$ is its field of fractions. 

In this case arguments used in the last part of the
previous proof  don't work because our functor is not defined
on $k$-algebras.  Moreover, the objects $\f(S)$ and
$\fr(S)$ are not isomorphic in general. 

To avoid this problem we will assume that the functors $\f$ and $\fr$
satisfy the additional axiom: 

\itemitem{\bf E.}{\bf (extension property)}
Given an $R$-algebra $R  \xra{i} S$ essentially smooth over the field $k$, 
given an $A$-algebra $i_S: A \ra S$ and 
an augmentation $\e: S \ra R$ of $i$, 
such that the following diagram commutes:
$$
\CD
A @>{i_S}>> S \\
@| @VV{\e}V \\
A @>>{i_R}> R
\endCD
$$

and given a multiplicative system $M$ with respect to a finite set
$\{\m_i\}_{i\in \I}$ of maximal ideals in $S$ 
with the property $\e^{-1}(\m_R)\subset \cup_{i\in \I}\m_i$

there exist

\itemitem{a)} the localization $S_g$ for a certain $g\in M$
with a finite etale extension $e: S_g \ra \St$;
\itemitem{b)} an augmentation $\et: \St \ra R$ for the
inclusion $R \xra{i} S_g \xra{e} \St$ such that the following
diagram commutes: 
$$
\CD
S_g @>{e}>> \St \\
@V{\e}VV @VV{\et}V \\
R @= R
\endCD
$$
\itemitem{c)} a natural transformation $\Phi: \f \ra \fr$
between two functors $\f$ and $\fr$ restricted to the
category of $\St$-algebras along $A \xra{i_S} S_g \xra{e} \St$ and   
$R \xra{i} S_g \xra{e} \St$
correspondingly, such that the morphism
$
\Phi(\St \xra{\et}R): \f(\St \xra{\et} R) \ar \fr(\St \xra{\et}R)
$
is the identity.

In particular, the following diagram commutes:
$$
\CD
\f(\St) @>{\Phi(\St)}>> \fr(\St) \\
@V{\et^*}VV   @VV{\et^*}V \\
\f(R) @= \fr(R)
\endCD
$$
where $R$ is considered as the $\St$-algebra via the augmentation $\et$.

\subhead{\bf Remark~3}\endsubhead
Since $\e^{-1}(\m_R)\cap M=\varnothing $, 
we can extend our augmentation $\e$ given on $S$ 
to the augmentation given on the localization $S_g$, 
for any $g \in M$, i.e., we have the commutative diagram:
$$
\CD
S @>{i_g}>> S_g \\
@V{\e}VV @VV{\e}V \\
R @= R
\endCD
$$
with the canonical inclusion $S \hra{i_g} S_g$.
Indeed, since $R$ is local with the unique maximal ideal $\m_R$, 
the image $\e(g)$ of any $g\in M$ is invertible in $R$. \qed

Then our main result can be stated as follows:

\proclaim{Theorem {\rm (Non-Constant Case)}}
Let $\f$ be the functor on the category of $A$-algebras and 
let $\fr$ be its restriction
to the category of $R$-algebras.
If $\f$ and $\fr$ satisfy axioms {\bf C}, {\bf TE},
{\bf TA}, {\bf TB} and {\bf E}
then the homomorphism $\f(R) \ra \f(K)$ 
induced by the canonical inclusion is injective.
\endproclaim

As in the proof of the Constant case (see the beginning of subsection~1.3),
by using continuity {\bf C}, we can reduce the problem 
about the injectivity of the map $\f(R)\ra \f(K)$ to the following one:
\proclaim{Proposition}
For a given  $\a \in \f(A)$ and $f \in \p$ 
such that the image $\a_f$ in $\f(A_f)$ of $\a$ is trivial
the element $i_R^*(\a)$ is trivial in $\f(R)$.
\endproclaim

\demo{Proof} As in Constant Case the proof consists of two steps.
On the first step for any triple $(A,R,f)$ we build up a $3 \times 5$ commutative diagram.
On the second step by playing with this diagram for the specially chosen triple
we complete the proof of the proposition.

In contrast with Constant Case the first step is more complicated and 
can be subdivided as follows:

First we produce starting data to apply axiom {\bf E}.
In particular we construct an {$R$-algebra} $S$ and a multiplicative system $M$.
Second we construct a lot of elements $h \in M$ such that the extension $S_h/R$
satisfies properties {\bf S1}, {\bf S2} and {\bf S3} of subsection~1.1.
Third for the specially chosen $h\in M$ we construct an extension $\St_h/R$ satisfying
conditions {\bf S1}, {\bf S2}, {\bf S3} and such that the diagram \ssp{} below commutes.
And finally we build up the  $3 \times 5$ commutative diagram.
\enddemo 

\subhead{I}\endsubhead
Let $R \xra{i} S$ be the extension built up by using
Quillen's Trick for any triple $(A,R,f)$. 
In particular, it means that there is a commutative diagram
$$
\CD
A @>{i_S}>> S \\
@| @VV{\e}V \\
A @>>{i_R}> R
\endCD
\tag{$*$}
$$
where $\e$ is the augmentation from {\bf S2}, $S$ is finite over $R[t]$ 
for specially chosen $t$ in $S$ and 
$S/(f)$ is finite over $R$ by {\bf S1}.

Let $I=\ker(\e)$, $J=(f)\cap I$ and $J'=(J\cap R[t])S$ 
(we identify $R[t]$ with its embedding in $S$).
Observe, that $J'\subset J \subset I$. 

\proclaim{Lemma 1} 
The quotients $S/J$ and $S/J'$ are finite over $R$.
\endproclaim

\demo{Proof}
Since $J=(f) \cap I$, there is the inclusion 
$S/J \hookrightarrow S/(f) \times S/I$ induced by projections.
It is easy to see now that  $S/J$ is finite over $R$.

Denote $J_t= J \cap R[t]$ the ideal in $R[t]$. 
Then $S/J'=S/J_tS=S\otimes_{R[t]}(R[t]/J_t)$ and,
since $S$ is finite over $R[t]$, 
our quotient $S/J_tS$ is finite over $R[t]/J_t$.
But there is the inclusion $R[t]/J_t \hookrightarrow S/J$ 
induced by the given one $R[t] \hookrightarrow S$ and
we have just proved that $S/J$ is finite over $R$. 
Hence, $R[t]/J_t$ is finite over $R$ and we have done. \qed
\enddemo

Let $\{\m_i\}_{i\in \I}$ be the system of all maximal
ideals lying over $J'$.  It is finite because of the lemma.
Let $M$ be the corresponding multiplicative system.
We see that $\e^{-1}(\m_R)$ lies above $I\supset J'$, 
so it is one of the $\m_i$, 
and $\e^{-1}(\m_R) \subset \cup_{i\in \I}\m_i$. 

\subhead{II}\endsubhead
Now we are interesting in the properties of localizations 
taken via the multiplicative system~$M$.

\proclaim{Lemma 2}
For any localization $S_g$, $g \in M$,
{\rm i)}~the quotient $S_g/(f)$ is finite over $R$; \ \  
{\rm ii)}~there exists an extension  $S_h$ of $S_g$, $h\in M$, 
such that $S_h$ has a structure of a finite algebra over $R[t']$ 
for some specially chosen~$t' \in S_h$.
\endproclaim

\demo{Proof} 
{\rm i)}~Since $(f) \supset J$, it is enough to show that $S_g/JS_g$
is finite over $R$. By construction of $M$ we know that $g
\in M$ is invertible in $S/J$, therefore,
$S_g/JS_g=(S/J)_g=S/J$ is finite over $R$.  

{\rm ii)}~Indeed, it is the reformulation of Lemma 8.2 \cite{PO2}.  \qed
\enddemo

Now show that the extension $S_h/R$ satisfies properties 
{\bf S1}, {\bf S2} and {\bf S3} of subsection~1.1. 

The property {\bf S1} follows from  i) and ii). 
The existance of augmentation {\bf S2} can be conclude from Remark~3. 
The property {\bf S3} is the consequence of the fact that 
the localization of a smooth algebra is again smooth and 
the fact that the element $h$ is coprime with the kernel of
augmentation~$I$. 

\subhead{III}\endsubhead 
Apply now axiom {\bf E} to the data $(R \xra{i} S, A \xra{i_S} S, S \xra{\e} R, M)$.
We get that there exists the localization $S_g$, $g\in M$, 
the finite etale extension $e: S_g \ra \St$, the augmentation
$\et:\St \ra R$  and the natural transformation $\Phi: \f \ra \fr$ 
between two functors $\f$ and $\fr$
restricted to the category of $\St$-algebras satisfying properties {\bf E.(b)} and {\bf E.(c)}.

By the previous step, we can find an extension $S_h$, $h\in M$, of $S_g$
with the property that $S_h/R$ satisfies conditions {\bf S1}, {\bf S2}, {\bf S3} of subsection~1.1.

Consider now the canonical tensor product diagram:
$$
\CD
S_h @>{i_1}>> S_h\otimes_{S_g} \St \\
@AAA  @AA{i_2}A \\
S_g @>>{e}> \St
\endCD
$$
where $i_1: s_h \mapsto s_h\otimes 1$ and 
$i_2: \tilde{s} \mapsto 1\otimes \tilde{s}$ 
are the canonical inclusions.

Set $\St_h=S_h\otimes_{S_g} \St$. 
We claim that the extension $\St_h/R$ satisfies 
properties {\bf S1}, {\bf S2}, {\bf S3}.

Indeed, we already know that the extension $S_h/R$ does satisfy. 
Define the augmentation map as 
$\et_h: s_h\otimes \tilde{s} \mapsto \e(s_h)\et(\tilde{s})$, thus,
we have checked {\bf S2}. 
Since $\St_h/S_h$ is the finite etale, the properties {\bf S1} and {\bf S3} hold.

Clearly, we have the commutative diagram:
$$
\CD
\St @>{i_2}>> \St_h \\
@V{\et}VV @VV{\et_h}V \\
R @= R
\endCD
\tag{$*'$}
$$

In contrast with the extension $S_h/R$ for the extension $\St_h/R$
by {\bf E.(c)} we have the commutative diagram
$$
\CD
\f(\St_h) @>{\Phi(\St_h)}>> \fr(\St_h) \\
@V{\et_h^*}VV  @VV{\et_h^*}V \\
\f(R) @= \fr(R)
\endCD
\tag{$**'$}
$$

\subhead{\bf IV}\endsubhead 
Consider the modification of the diagram \sss{} built
up by using steps {\bf I}--{\bf III} for any triple $(A,R,f)$:
$$
\CD
\f(A_f) @>>> \f(S_{gf}) @>>> \f(\St_f) @>>> \f(\St_{hf})
@>{\Phi(\St_{hf})}>> \fr(\St_{hf}) \\ 
@AAA  @AAA  @AAA  @AAA  @AAA \\
\f(A) @>{i_S^*\circ i_g^*}>> \f(S_g) @>{e^*}>> \f(\St)
@>{i_2^*}>> \f(\St_h) @>{\Phi(\St_h)}>> \fr(\St_h) \\ 
@|  @V{\e^*}VV  @V{\et^*}VV  @V{\et_h^*}VV  @VV{\et_h^*}V \\
\f(A) @>>{i_R^*}> \f(R) @= \f(R) @= \f(R) @= \fr(R)
\endCD
$$
where the lower squares are constructed as follows:
the left side is induced by ($*$) and
by the commutative diagram from Remark~3; the next squares are induced by {\bf E.(b)},
($*'$) and ($**'$), correspondingly. 
Hence, we get the desired $3 \times 5$ commutative diagram.

\subhead{\bf V}\endsubhead
Repeat now the arguments used
in the step {\bf II} of the proof of Constant Case (see the end of subsection~1.3): 

Let $\a \in \f(A)$ and $f$ be an element from $\p$ such
that $\a_f=0$. 

Denote by $\a_{\St_h}$ the image in $\fr(\St_h)$ of $\a$ under the composition 
$i_S^*\circ i_g^*\circ e^*\circ i_2^* \circ \Phi(\St_h)$. 
Since the image $\a_f$ in $\f(A_f)$ of $\a$ is trivial, then
by the commutativity of the diagram 
the image in $\fr(\St_{hf})$ of the element $\a_{\St_h}$ is trivial, too.

Hence, we can apply Specialization Lemma to the right
column of our diagram and  we get ${\et_h}^*(\a_{\St_h})=0$. 
By the commutativity of the diagram $i_R^*(\a)={\et_h}^*(\a_{\St_h})$, thus,  $i_R^*(\a)=0$.
And we have finished. \qed

% ***************************************************
% ***************************************************
% ***************************************************
\head{3. Applications}\endhead

We keep all notations and agreements used before. 
Let $A$ be a smooth algebra over an infinite field $k$ of 
characteristic different from 2 and let
$R$ be its local regular ring at some prime ideal $\p$. 
Let $K$ be a field of fractions of $R$. 

In the present section we will apply the theorem (Non-Constant Case)
to some specially chosen functor $\f$ given on  the category of $A$-algebras.

Namely, the point is to show that this functor satisfies axioms 
{\bf C}, {\bf TE}, {\bf TA}, {\bf TB} and {\bf E} from the previous sections.
Then we can use our  main theorem (see section~2) 
\proclaim{Theorem {\rm (Non-Constant Case)}}
Let $\f$ be the functor on the category of $A$-algebras and 
let $\fr$ be its restriction
to the category of $R$-algebras.
If $\f$ and $\fr$ satisfy axioms {\bf C}, {\bf TE}, {\bf TA}, {\bf TB} 
and {\bf E}
then the homomorphism $\f(R) \ra \f(K)$ 
induced by the canonical inclusion is injective.
\endproclaim
Hence, we get the injectivity of the map $\f(R) \ra \f(K)$ 
(see the (a) statement of the theorems below).
Finally, playing with a few long exact cohomology sequences and
using the injectivity above 
we get the injectivity on the first cohomology level 
(see the (b) statement of the theorems below).

As a consequence, 
we will get two cases of Grothendieck's Conjecture about
principal homogeneous spaces \cite{Gr} ---  
the case of Special Linear group and the case of Special Unitary group.  

Now let describe our special functors and results we are interested in
(the detailed proofs one will find in the corresponding subsections below):

\subsubhead{Linear Case}\endsubsubhead
Let $\A$ be some Azumaya algebra over the given local regular ring $R$ (for
the definition see subsection 3.1). 
Let $\nrd: \A^* \ra R^*$ denotes the reduced norm homomorphism.
For any $R$-algebra $T$ let $\A_T=\A\otimes_R T$ be the extended Azumaya algebra over $T$.

Define the group scheme $\SL_{1,\A}$ related to the Azumaya $R$-algebra $\A$ as
$$
\SL_{1,\A}: T \mapsto \SL(\A_T)=\{ a \in \A_T \mid  \nrd(a)=1 \}.
$$

Our aim is to show: 

\proclaim{Theorem }
Let $\A$ be an Azumaya algebra over a local regular ring $R$ of geometric type. Then 

\itemitem{\rm (a)} the homomorphism  $R^*/\nrd(\A^*) \ar K^*/\nrd(\A_K^*)$ is injective;

\itemitem{\rm (b)} the canonical map $\hml(R, \SL_{1,\A}) \ar \hml(K, \SL_{1, \A_K})$
on the first cohomology groups induced by the canonical inclusion is injective.
\endproclaim
\demo{Proof} See subsection 3.2 \enddemo

Observe that the statement (a) of the theorem corresponds to the case when the functor
$\f$ is defined as $\f: T \mapsto T^*/\nrd(\A_T^*)$.

We recall that the conjecture of Grothendieck \cite{Gr} states the triviality of the kernel
of the canonical map
$$
\hml(R,G) \ar \hml(K,G_K)
$$
for any reductive flat group scheme $G$ over a regular semilocal ring $R$. 
Thus, the statement
(b) of the theorem above proves the following assertion:

\proclaim {Corollary {\rm (Special Linear Case)}}
The Grothendieck's conjecture is true for the group scheme $G=\SL_{1,\A}$ related
to an Azumaya algebra $\A$ over a local regular ring $R$ of geometric type.
\endproclaim

In the same notation assume that there is the additional structure on our Azumaya algebra:
Let $(\A,\s)$ be an Azumaya algebra with involution over $R$ (for the definition see subsection 3.1).
 
We know that there can be three different types of involution: ortogonal, sympletic and unitary.
It turns out that in the case of ortogonal and sympletic involution we can prove the same 
results as above by using well-known facts about quadratic forms 
(for the details see subsection 3.3).  
So the only interesting case for us is the unitary case:

\subsubhead{Unitary Case}\endsubsubhead
Let $(\A,\s)$ be an Azumaya algebra with unitary involution over $R$.
It means that there is a tower $\A/C/R$, where $C$ is the center of $\A$ and 
$C/R$ is an etale quadratic extension over $R$ with restricted involution $\s$.
Therefore, $C_K/K$ is a separable quadratic extension of the corresponding fields of fractions.

Let $\U(\A_T)=\{a \in \A_T \mid aa^\s=1 \}$ be the unitary group of an algebra $(\A_T,\s)$ for
any $R$-algebra $T$.
We  define the group scheme $\SU_{1,\A}$ related to the Azumaya $R$-algebra $\A$ with unitary involution 
$\s$ as
$$
\SU_{1,\A}: T \mapsto \SU(\A_T)=\{a \in \A_T \mid aa^\s=1,\; \nrd(a)=1\}
$$

Our goal will  be to show:
\proclaim{Theorem } 
Let $(\A,\s)$ be an Azumaya algebra with unitary involution
over a local regular ring $R$ of geometric type. 
Then 

\itemitem{\rm (a)} the homomorphism  $\U(C)/\nrd(\U(\A)) \ar \U(C_K)/\nrd(\U(\A_K))$ is injective;

\itemitem{\rm (b)} the kernel of the canonical map $\hml(R, \SU_{1,\A}) \ar \hml(K, \SU_{1, \A_K})$ is trivial.
\endproclaim
\demo{Proof} See subsection 3.4 \enddemo

Clearly, the statement (a) respects the functor $\f: T \mapsto \U(C_T)/\nrd(\U(\A_T))$, 
where $C_T=C\otimes_R T$ is the center of the extended Azumaya algebra $\A_T$.

Thus, the statement
(b) of the theorem above proves the following assertion:

\proclaim {Corollary {\rm (Special Unitary Case)}}
The Grothendieck's conjecture is true for the group scheme $G=\SU_{1,\A}$ related
to an Azumaya $R$-algebra $\A$ with unitary involution 
over a local regular ring of geometric type.
\endproclaim

\subsubhead{Torsion Cases}\endsubsubhead
Next two theorems represent independent interest. One can look on it
as on the modification of the corresponding Linear and Unitary Case.

For the functor $\f: T \mapsto T^*/\nrd(\A_T^*)(T^*)^d$, $d \in \N$,
we have
\proclaim{Theorem {\rm (Linear Torsion Case)}}
Let $\A$ be an Azumaya algebra over a local regular ring $R$ of geometric type 
and $d$ be some natural number.
Then the homomorphism 
$$
R^*/\nrd(\A^*)(R^*)^d \ar K^*/\nrd(\A_K^*)(K^*)^d
$$
is injective.
\endproclaim

And for the functor $\f: T \mapsto \U(C_T)/\nrd(\U(\A_T))\U(C_T)^d$, $d \in \N$, we have
\proclaim{Theorem {\rm (Unitary Torsion Case)}}
Let $(\A,\s)$ be an Azumaya algebra with unitary involution
over a local regular ring $R$ of geometric type and $d$ be some natural number. 
Then the homomorphism  
$$
\U(C)/\nrd(\U(\A))\U(C)^d \ar \U(C_K)/\nrd(\U(\A_K))\U(C_K)^d
$$ 
is injective.
\endproclaim
\demo{Proof} 
The torsion cases  one can prove easily by following
the proof of the corresponding Linear and Unitary case (see subsections 3.2 and 3.4). 
 \enddemo

\subhead\nofrills{3.1. The Properties of Azumaya algebras} \endsubhead

Let us formulate some properties of Azumaya algebras which we are going to use.
Mostly, one can find them in   \cite{Kn} and \cite{PO2}.

We recall that an Azumaya $R$-algebra over a ring $R$ is an $R$-algebra $\A$ which
satisfies the following two properties: it is finitely generated projective $R$-module, and
the canonical $R$-algebra homomorphism $\A\otimes_R \A^{op} \ra \End_R(\A)$ is an isomorphism.
Clearly, for any $R$-algebra $T$, the scalar extension $\A_T=\A\otimes_R T$ of the algebra $\A$
is an Azumaya algebra over $T$.

We will use the reduced norm homomorphism $\nrd_R: \A^* \ra R^*$. 
This homomorphism respects the scalar extensions. 
In the case where $R$ is a field and $T$ is its algebraic closure, 
the reduced norm homomorphism $\nrd_T$ can be identified 
with the usual determinant $\det: \GL_d(T) \ra T^*$; 
this identification is induced by an $T$-algebra isomorphism 
between $\A_T$ and  the matrix algebra $M_d(T)$ mentioned above.

By an Azumaya algebra with involution over a ring $R$ we mean a pair
$(\A,\s)$ consisting of an $R$-algebra $\A$ and an $R$-linear involution $\s$ on $\A$ such that:
i)~$\A$ is an Azumaya algebra over its center $C$; \ \ 
ii)~$C$ is either $R$ or an etale quadratic extension of $R$; \ and\  
iii)~$C^\s=R$.
Observe that the involution $\s$ commutes with the reduced norm map, 
i.e. $\nrd_C(a^\s)=\nrd_C(a)^\s$ for any $a \in \A^*$.

When $C=R$ we have involution of the first kind. 
We say that involution $\s$ is of ortogonal (or sympletic) type if the dimension over $R$
of the set of symmetric elements  $\{a \in \A \mid a^\s=a \}$ of the algebra $\A$
equals $\tfrac{n(n+1)}{2}$ (or $\tfrac{n(n-1)}{2}$),
where $n$ is the degree of the Azumaya algebra $\A$ (see \cite{Kn}).

In the case when $C$ is quadratic etale over $R$ we have unitary involution
(involution of the second kind).

Now we fix a local regular ring $R=A_\p$ of a smooth $k$-algebra $A$. 
Then we have the following properties:
 
\subsubhead{\bf A1}\endsubsubhead
Since an Azumaya $R$-algebra $\A$ is given by the finite number of 
generators and relations we can find
the localization $A_g$, $g\in A\setminus \p$, such that the algebra $\A$ come from some Azumaya
$A_g$-algebra $\A_g$, i.e. $\A=\A_g\otimes_{A_g} R$. On geometric language it means
that we can extend an Azumaya algebra given at a point  to some neighbourhood of this point.

By the same reasons as before we can state that if 
there is an isomorphism $\Psi:\A \ra \B$ of Azumaya $R$-algebras
then there exists the localization $A_g$, $g \in A \setminus \p$, 
the Azumaya $A_g$-algebras $\A_g$ and $\B_g$, where
$\A=\A_g\otimes_{A_g} R$ and $\B=\B_g\otimes_{A_g} R$,
and the isomorphism $\Psi_g: \A_g \ra \B_g$ of the Azumaya $A_g$-algebras such that
$\Psi=\Psi_g \otimes_{A_g} \id_R$.

Observe that the arguments above also work in the case 
of Azumaya algebras with involutions and
when $R$ is a semilocal regular ring, i.e.
$R$ is the localization at some finite number of prime ideals of $A$.

\subsubhead{\bf A2}\endsubsubhead
We also will need in the following reformulation of the Proposition 7.1 of \cite{PO2}:
Let $\O$ be some semilocal regular ring such that there is the inclusion $R \hra{i} \O$ with
the augmentation $\e$.
Let $(\A_{\O},\s)$ and $(\B_{\O},\tau)$ be two Azumaya algebras with involution over $\O$, of the same rank.
Assume that there exists  an isomorphism $\psi: (\A_R,\s) \ra (\B_R,\tau)$. Then there exists
a finite etale extension $e: \O \ra \Ot$, an augmentation $\et: \Ot \ra R$ of $e$ over $R$ and an
isomorphism $\Psi: (\A_{\Ot},\st) \ra (\B_{\Ot},\tilde{\tau})$ 
of extended Azumaya algebras with involutions over $\Ot$ 
such that the extension 
$\Psi\otimes_{\Ot} \id_R: (\A_{\Ot},\st)\otimes_{\Ot} R \ra (\B_{\Ot},\tilde{\tau})\otimes_{\Ot} R$ 
of $\Psi$ via $\et$  coinsides with $\psi$.

%------------------------------------------------------------------------------------------------------------------
\subhead\nofrills{3.2. The Proof of Linear Case}\endsubhead

Let $R=A_\p$  be the local regular ring of the smooth $k$-algebra $A$ and 
let $\A$ be the Azumaya algebra over $R$.

Let $\f$ be the functor defined on the category of $R$-algebras as:
$$
\f: T \mapsto T^*/\nrd(\A_T^*).
$$

First of all, by property {\bf A1}, we may assume that 
our Azumaya algebra $\A$ over the ring $R$ come from 
some Azumaya algebra over the localization $A_g$, 
for some $g \in A \setminus \p$. 
Thus, our functor is defined on the category of $A_g$-algebras.

Since $A_g$ is smooth $k$-algebra and $R$ is again its local regular ring, 
we may write $A$ instead of $A_g$ and nothing will be changed. 
Hence, we will write $\A_R$ instead of $\A$ meaning that 
$\A_R$ is the scalar extension of the Azumaya $A$-algebra $\A$ 
via the canonical inclusion $A \hra{i_R} R$.
Thus, we may assume that our functor $\f$ is defined on the category of $A$-algebras.

Show that the functor $\f$ and its restriction $\fr$ to the category of
$R$-algebras satisfy axioms {\bf C}, {\bf TE}, {\bf TA}, {\bf TB} and {\bf E}:

\subhead{C}\endsubhead
Observe that the functor $G_m: T \mapsto T^*$ sending
any $A$-algebra $T$ to its group of units is continuous.
Since $\A$ is finitely generated projective as the $A$-module, 
the functor of extension of scalars $T \mapsto \A_T$ is continuous.
Thus, the functor $T \mapsto \nrd(G_m(\A_T))$ is continuous, too, and
we get the required.

\subhead{E}\endsubhead
Let we are under the hypothesis of the axiom {\bf E} (see section 2):
Let we have the $R$-algebra $R \xra{i} S$, 
the inclusion $i_S: A \ra S$,
the augmentation $\e: S \ra R$ of $i$
and the multiplicative system $M$.

For the functors $\f$ and $\fr$ restricted to the category of $S$-algebras 
we may write: 
$$
\f: T \mapsto T^*/\nrd(\A_T^*) \qquad\text{and}\qquad \fr: T \mapsto T^*/\nrd(\B_T^*),
$$
where for any $S$-algebra $T$,
$\A_T$ is the extension of scalars of $\A$ via the inclusion $A \xra{i_S} S \ra T$ and
$\B_T$ is the extension of scalars of $\A$ via $A \xra{i_R} R \xra{i} S \ra T$
(in general, $\A_T$ and $\B_T$ are not isomorphic).

We will check axiom {\bf E} in four steps:

First, by using  {\bf A2} property we will produce 
the finite etale extension $\O \xra{e} \Ot$ of the localization $\O=M^{-1}S$ and
the augmentation $\et: \Ot \ra R$,
such that the extended Azumaya $\Ot$-algebras $\A_{\Ot}$ and $\B_{\Ot}$ 
become equivalent via the isomorphism $\Psi$.

Secondly, we will show that the semilocal ring $\Ot$ 
is, indeed, the localization $M^{-1}\St$ of some $R$-algebra $\St$ 
which is finite etale over $S$.

After that, by using {\bf A1} we will find the localization $\St_g$, $g\in M$, of $\St$
such that the extended Azumaya $\St_g$-algebras $\A_{\St_g}$ and $\B_{\St_g}$
are still equivalent via the isomorphism $\Psi_g$ with $\Psi=\Psi_g\otimes_{\St_g}\Ot$. 
Considering the extension $e: S_g \ra \St_g$ and the augmentation 
$\et: \St_g \ra \Ot \xra{\et} R$ we get conditions (a) and (b) of axiom {\bf E}.

We will end with showing condition (c) of axiom {\bf E}.
For this purpose by using the isomorphism $\Psi_g: \A_{\St_g} \ra \B_{\St_g}$
we construct the natural equivalence $\Phi$ of the functors $\f$ and $\fr$
restricted to the category of $\St_g$-algebras. 
It turns out that all necessary properties for condition (c) are satisfied.

\subhead{\bf I}\endsubhead
We introduce the semilocal ring $\O$ as the localization $M^{-1}S$ of $S$. 
By  Remark~3 we have the
augmentation map $\e : \O \ra R$ for the inclusion $R \xra{i} S \hra{i_\O} \O$ 
compatible with augmentation on $S$.
Hence, there is the commutative diagram:
$$
\CD
A @>{i_S}>> S @>{i_\O}>> \O \\
@| @V{\e}VV @VV{\e}V \\
A @>>{i_R}> R @= R
\endCD
$$
Consider the extended Azumaya $\O$-algebras $\A_{\O}$ and $\B_{\O}$. 
Since the extensions $\A_R$ and $\B_R$ of the algebras 
$\A_{\O}$ and $\B_{\O}$ via the augmentation $\e$ coinside (see the diagram above),
we are under the hypothesis of property {\bf A2} (our $\psi$ is the identity).

So there exists
the finite etale extension $e: \O \ra \Ot$, 
the augmentation $\et: \Ot \ra R$ for the inclusion 
$R \xra{i} S \xra{i_\O} \O \xra{e} \Ot$ compatible with $\e$
and the isomorphism $\Psi: \A_{\Ot} \ra \B_{\Ot}$ of extended Azumaya algebras 
such that its extension $\Psi \otimes_{\Ot} \id_R$ via $\et$ is the identity.

\subhead{\bf II}\endsubhead
By the properties of finite etale extensions
there exists the localization $S_h$, $h\in M$, a finite etale extension $e: S_h \ra \St$
such that $\Ot=\O\otimes_{S_h}\St$. Since $\O=M^{-1}S_h$, we have $\Ot=M^{-1}\St$.
To simplify the notation we will write $S$ instead of $S_h$. Thus, we get the commutative
diagram:
$$
\CD
\St @>localization>via \; M> \Ot \\
@A{e}AA @AA{e}A \\
S @>localization>via \; M> \O
\endCD
$$

\subhead{\bf III}\endsubhead
By the properties of finite extensions
for the direct system of localizations $\{S_g\}_{g\in M}$ with $\varinjlim_{g\in M}S_g = \O$ 
we have an induced direct system $\{\St_g\}_{g\in M}$, 
where $\St_g=S_g\otimes_S\St$, and
the canonical map $\varinjlim_{g\in M} \St_g \ra \Ot$ is an isomorphism.
Thus, one has the diagram
$$
\CD
\St @>>> \St_g @>>> \Ot \\
@A{e}AA  @A{e}AA  @AA{e}A \\
S @>>> S_g @>>> \O
\endCD
$$
where $e$ is finite etale.

Since our Azumaya $\Ot$-algebras $\A_{\Ot}$ and $\B_{\Ot}$ are isomorphic via $\Psi$, 
applying {\bf A1} to the case $A=\St$ and $R=\Ot$ we get that
there exists a localization $\St_g$, $g\in M$, such that the 
Azumaya $\St_g$-algebras $\A_{\St_g}$ and $\B_{\St_g}$ are isomorphic via $\Psi_g$.

Take the composition $\et:\St_g \ra \Ot \xra{\et} R$ to be the augmentation on $\St_g$.
Clearly, it is compatible with $\e$. Since the isomorphism $\Psi$ is got from $\Psi_g$ 
by extension of  scalars $\St_g \ra \Ot$  and the extension of $\Psi$ 
via $\et: \Ot \ra R$ is the identity, the extension of $\Psi_g$ 
via the augmentation $\et: \St_g \ra R$  is the identity as well.

To simplify the notations we replace $\St_g$ by $\St$ and $\Psi_g$ by $\Psi$.
So we have constructed:
\itemitem{a)} the localization $S_g$, $g\in M$, and the finite etale extension $S_g \xra{e} \St$;
\itemitem{b)} the augmentation $\et$ on $\St$ compatible with $\e$.

Thus, we get conditions (a) and (b) of axiom {\bf E}.

\subhead{\bf IV}\endsubhead
The last step is to check {\bf E.(c)}:

Since we have the isomorphism $\Psi: \A_{\St} \ra \B_{\St}$ on Azumaya $\St$-algebras, 
we have the isomorphism $\Psi(T): \A_T \ra \B_T$ of Azumaya $T$-algebras 
got by extension of scalars for any $\St$-algebra $T$.
By definition there is the commutative diagram:
$$
\CD
\A_T^* @>\Psi(T)>> \B_T^* \\
@V{\nrd}VV @VV{\nrd}V \\
T^*=Z(\A_T)^* @>>\Psi_Z(T)> Z(\B_T)^*=T^*
\endCD
$$
where $\Psi_Z(T)$ is the restriction of $\Psi(T)$ 
to the groups of units $T^*$ of the centers  
of the Azumaya algebras $\A_T$ and $\B_T$. 
Thus, there is the isomorphism on the quotients
$$
\Phi(T)=\overline{\Psi_Z(T)}: T^*/\nrd(\A_T^*) \ra T^*/\nrd(\B_T^*).
$$
and  $\Phi: \f \ra \fr$ is the natural equivalence of the functors $\f$ and $\fr$
on the category of $\St$-algebras.

Moreover, the extension of the isomorphism $\Phi(\St): \f(\St) \ra \fr(\St)$
via the augmentation $\et: \St \ra R$ is the identity,
i.e. $\Phi(R)=\id_{\f(R)}: \f(R) \ra \fr(R)$. 

Observe that in our concrete case the isomorphism $\Psi_Z(T)$ is the identity. 
In general, it is not true. 
For instance,  when we have Azumaya algebras with unitary involutions
(see subsections 3.1 and 3.4)
the isomorphism on the centers may not coinside with the identity.

\subhead{TE}\endsubhead
Let $T$ be an $R$-algebra finitely generated projective as the $R$-module. 
Thus, $T$ is a semilocal ring.
We define the transfer map 
$$
\tr_R^T: T^*/\nrd(\A_T^*) \ar R^*/\nrd(\A_R^*)
$$ 
to be the usual norm map $\n_R^T :T^* \ra R^*$ on the quotients.

To see that it is well-defined we need in:
\proclaim{Lemma 3}
Let $\A$ be an Azumaya algebra over a semilocal ring $R$. 
Let  $T$ be an $R$-algebra finitely generated projective as the $R$-module and
let  $\n_R^T: T^* \ra R^*$ be its norm map. 
Then we claim that:
$$
\n_R^T(\nrd_T(\A_T^*)) \subset \nrd_R(\A^*).
$$
\endproclaim

\demo{Proof}To prove this inclusion we are passing our reduced norms through $K_1$-groups.
\enddemo
By results of \cite{Ba} (ch.V, Theorem 9.1) there is the commutative diagram:
$$
\CD
\A_T^* @>{i_T}>> K_1(\A_T) @>{\n^*}>> K_1(\A) @<{i}<{surj}< \A^* \\
@V{\nrd_T}VV  @V{\nrd_T^*}VV  @VV{\nrd_R^*}V  @VV{\nrd_R}V \\
T^*  @=  T^*  @>>{\n}>  R^*  @=  R^*
\endCD
$$
where $\n^*$ is the norm map on $K_1$-groups; $i^*$ is the surjective 
homomorphism induced by inclusion $\A^* = \GL_1(\A) \hra{} K_1(\A)$; 
and $\nrd^*$ is the reduced norm
on $K_1$ of an Azumaya algebra. If $\A$ splits, i.e. $\A=M_d(R)$, we can write our reduced norm
as the composition:
$$
\CD
\nrd^*: K_1(M_d(R)) @>{\text{Morita equivalence}}>>  K_1(R) @>{\det}>>  R^*.
\endCD
$$

We get the required inclusion since $i$ is surjective. \qed

\subhead{TA}\endsubhead
The additivity follows from the corresponding property of the norm map:
If $T_1$ and $T_2$ are finitely generated projective over $R$ then
for any $(t_1,t_2)\in T_1\times T_2$ we have 
$\n_R^{T_1\times T_2}(t_1,t_2)=\n_R^{T_1}(t_1)\n_R^{T_2}(t_2)$.
Since $\fr(T_1\times T_2)=\fr(T_1) \times \fr(T_2)$, we get the required.

\subhead{TB}\endsubhead
Let $S$ be any $R[t]$-algebra finitely generated projective as the $R[t]$-module.
Observe that the functor $G_m: T\mapsto T^*$ with the usual norm in the
role of the transfer map satisfies homotopy invariance {\bf H} and 
base changing {\bf TB'} properties, thus, it satisfies {\bf TB} (see Remark~1).

To see the axiom {\bf TB} look on the  diagram induced
by {\bf TB} applying to $G_m$:
$$
\CD
G_m(S) @>{\can_0^*}>> G_m(S/(t)) \\
@V{\can_1^*}VV        @VV{\n_0}V \\
G_m(S/(t-1)) @>>{\n_1}> G_m(R)
\endCD\tag$\ast$
$$
where $\can_i^*$, $i=0,1$, are induced by the canonical projections, 
and $\n_i$ denote the corresponding norm map.

The following two diagrams induced by the canonical projection 
$$\pr: G_m(T)=T^* \ar T^*/\nrd(\A_T^*)=\fr(T)$$  
are commutative as well:
$$
\CD
G_m(S) @>{\can_i^*}>> G_m(S/(t-i)) @>{\n_i}>> G_m(R) \\
@V{\pr}VV @VV{\pr}V @VV{\pr}V \\
\fr(S) @>>{\can_i^*}> \fr(S/(t-i)) @>>{\tr_i}> \fr(R)
\endCD
$$
where $\tr_i$, $i=0,1$, denote the corresponding transfer map.

Thus, since the projection $\pr$ is surjective,
we get the required diagram commutes:
$$
\CD
\fr(S) @>{\can_0^*}>> \fr(S/(t)) \\
@V{\can_1^*}VV        @VV{\tr_0}V \\
\fr(S/(t-1)) @>>{\tr_1}> \fr(R)
\endCD
$$
And we have checked that the functors $\f$ and $\fr$ satisfy axioms 
{\bf C},  {\bf TE}, {\bf TA}, {\bf TB}, {\bf E}. \qed

Applying now our theorem (Non-Constant Case) to the functors $\f$ and $\fr$ 
we get the statement (a) of Linear Case:
\proclaim{Theorem {\rm (a)}}
The map $R^*/\nrd(\A^*) \ar K^*/\nrd(\A_K^*)$ induced by the canonical inclusion is injective.
\endproclaim

To prove the statement (b) let look on the short exact
sequence of smooth group schemes:
$$
1 \ar \SL_{1,\A} \ar \GL_{1,\A} \overset{\nrd}\to{\ar} G_m \ra 1,
$$
where $\GL_{1,\A}: T \mapsto \A_T^*$ for any $R$-algebra $T$.

It induces the long exact cohomology sequence:
$$
1 \ar \SL(\A) \ar \A^* \overset{\nrd}\to{\ar} R^*
\ar \hml(R,\SL_{1,\A}) \ar \hml(R,\GL_{1,\A}) \ar \cdots
$$ 

Since $R$ is local, the group $\hml(R,\GL_{1,\A})$ is trivial (see \cite{Kn}).

The inclusion $R \hookrightarrow K$ induces the natural map
on our long cohomology sequence, so that the diagram
$$
\matrix
1 & \lra{} & \SL(\A) & \lra{} & \A^* & \lra{\nrd} & R^* & \lra{} & \hml(R,\SL_{1,\A}) & \lra{} & 1 \\
  &  & \da & & \da & & \da  & & \da & & \\
1 & \lra{} & \SL(\A_K) & \lra{} & \A_K^* & \!\!\!\! \underset{\nrd_K}\to{\longrightarrow}\!\!\!\! & K^* & \lra{} & \hml(K,\SL_{1,\A_K}) & \lra{} & 1 \\
\endmatrix
$$
commutes.

Taking the cokernels from the left side we get:
$$
\CD
R^*/\nrd(\A^*) @>{\cong}>> \hml(R,\SL_{1,\A}) \cr
@VVV @VVV \cr
 K^*/\nrd(\A_K^*) @>>{\cong}> \hml(K,\SL_{1,\A_K}) 
\endCD
$$

Since the left arrow is injective by (a) we get the statement (b):

\proclaim{Theorem {\rm (b)}}
The map $\hml(R,\SL_{1,\A}) \ar \hml(K,\SL_{1,\A_K})$ induced by the canonical inclusion is injective.
\endproclaim

\subhead\nofrills{3.3. Ortogonal and Sympletic Cases} \endsubhead

Let $(\A,\s)$ be an Azumaya algebra with ortogonal involution over $R$
(for the definition see subsection 3.1).

Let $\Or(\A_T)=\{a \in \A_T \mid aa^\s=1 \}$ be the ortogonal group 
of an algebra $(\A_T,\s)$ for any $R$-algebra $T$.
We  define the group scheme $\SO_{1,\A}$ related to the Azumaya $R$-algebra $\A$ 
with ortogonal involution $\s$ as
$$
\SO_{1,\A}: T \mapsto \SO(\A_T)=\{a \in \A_T \mid aa^\s=1,\; \nrd(a)=1\}.
$$

For the field of fractions $K$ we have $\Or(K)=\{x \in K \mid x^2=1\}=\{\pm 1\}$.
Since $R$ is the local ring of an affine variety over the field $k$,
we get also $\Or(R)=\{\pm 1\}$.

We would like to show the analogy of the statement (a) of Linear Case, i.e.
the injectivity of the map
$$
\Or(R)/\nrd(\Or(\A)) \ar \Or(K)/\nrd(\Or(\A_K))
$$
induced by the canonical inclusion.

First of all, note that our map is always surjective and it is not injective if and  only if
$\nrd(\Or(\A))=\{1\}$ and $\nrd(\Or(\A_K))=\{\pm 1\}$.

Let $\nrd(\Or(\A_K))=\{\pm 1\}$ then by theorem of Knezer \cite{BI} we conclude that $\A_K$ splits, 
i.e. $\A_K$ is the matrix algebra over $K$. 
Moreover, by theorem of Grothendieck \cite{BI}, the algebra $\A$ splits, too. 
Hence, we get $\nrd(\Or(\A))=\{\pm 1\}$ and our map must be injective anyway,
indeed, it is an isomorphism.

Now look on the short exact
sequence of smooth group schemes:
$$
1 \ar \SO_{1,\A} \ar \Or_{1,\A} \overset{\nrd}\to{\ar} \Or \ra 1.
$$

It induces the long exact cohomology sequence:
$$
1 \ar \SO(\A) \ar \Or(\A) \overset{\nrd}\to{\ar} \Or(R)
\ar \hml(R,\SO_{1,\A}) \ar \hml(R,\Or_{1,\A}) \ar \cdots
$$ 

The canonical inclusion $R \hookrightarrow K$ induces the natural map
on our long cohomology sequence, so that the diagram
$$
\matrix
1 & \lra{} & \SO(\A) & \lra{} & \Or(\A) & \lra{\nrd} & \Or(R) & \lra{} &  \hml(R,\SO_{1,\A}) & \lra{} & \hml(R,\Or_{1,\A}) \cr
 & & \da &  & \da & & \da & & \da && \da \cr
1 & \lra{} & \SO(\A_K) & \lra{} & \Or(\A_K) & \!\!\!\!\underset{\nrd_K}\to{\ar}\!\!\!\! & \Or(K) & \lra{} &  \hml(K,\SO_{1,\A_K}) & \lra{} & \hml(K,\Or_{1,\A_K}) 
\endmatrix
$$
commutes.

Taking the cokernels from the left side we get:
$$
\CD
1 @>>> \Or(R)/\nrd(\Or(\A)) @>>> \hml(R,\SO_{1,\A}) @>>> \hml(R,\Or_{1,\A})\\
@. @V{\cong}VV @VVV @VVV \\
1 @>>>  \Or(K)/\nrd(\Or(\A_K)) @>>> \hml(K,\SO_{1,\A_K}) @>>> \hml(K,\Or_{1,\A_K}) 
\endCD
$$

Since the left vertical arrow is the isomorphism
and the right vertical arrow has the trivial kernel by the result of \cite{PO2},
we get the triviality of the kernel of the middle vertical arrow.

Thus, we have proved the following assertion:

\proclaim {Theorem {\rm (Special Ortogonal Case)}}
The Grothendieck's conjecture is true for the group scheme $G=\SO_{1,\A}$ related
to an Azumaya algebra $\A$ with ortogonal involution over a local regular ring $R$ of geometric type.
\endproclaim

Special Sympletic Case of Grothendieck's conjecture about principal homogeneous
spaces follows immediately from \cite{PO2}.

\subhead\nofrills{3.4. The Proof of Unitary Case}\endsubhead

We keep all notations and definitions used in the subsection~3.1. 
The most arguments of our discussion here are taken from Linear Case. 
The only difference is that we have an additional structure of unitary involution
on our Azumaya algebras.

Thus, let $(\A,\s)$ be an Azumaya algebra with unitary involution over $R$.

Let $\f$ be the functor defined on the category of $R$-algebras as:
$$
\f: T \mapsto \U(C_T)/\nrd(\U(\A_T)),
$$
where $\U(C_T)=\{c\in C_T \mid cc^\s=1\}$ is the unitary group of the center of $\A_T$.

By the same arguments as in subsection 3.2, 
we may assume that our functor is given on the category of $A$-algebras.

Now the problem is to check the axioms {\bf C}, {\bf TE}, {\bf TA}, {\bf TB} and {\bf E}.
Axioms {\bf C} and {\bf E} can be proved following exactly to the
corresponding proofs in subsection~3.2.

The main difficulty is to show the existance of the transfer map for some 
finitely generated projective extension $T/R$.
Indeed, we would like to see the norm homomorphism in the role of the transfer map again but
there is no way to show the inclusion 
$$
\n_C^{C_T}(\nrd_{C_T}(\U(\A_T)))\subset \nrd_C(\U(\A))
$$
by using arguments with $K_1$ (there is no well-defined norm map for the unitary $K_1$).

The following important theorem gives us another possibility to do this.

\proclaim{Theorem {\rm (Norm Principle for the unitary group)}}
Let $T$ be a semilocal ring with infinite residue fields of characteristic different from 2.
Let $(\A_T,\s)$ be an Azumaya algebra with unitary involution $\s$ over $T$.
Let $C_T$ be the center of $\A_T$, so $C_T/T$ is an etale quadratic extension 
with the restricted involution $\s$. Then the following equality holds: 
$$
\nrd_{C_T}(\U(\A_T))=\nrd_{C_T}(\A_T^*)^{1-\s},
$$
where $c^{1-\s}=c(c^\s)^{-1}$, for any $c\in C_T^*$.
\endproclaim

\demo{Proof} See section 4 below.
\enddemo

Since the norm commutes with the involution we get that 
$$
N_C^{C_T}(\nrd_{C_T}(\U(\A_T)))=N_C^{C_T}(\nrd_{C_T}(\A_T^*)^{1-\s})=
$$
$$
N_C^{C_T}(\nrd_{C_T}(\A_T^*))^{1-\s}\subset \nrd_C(\A^*)^{1-\s}=\nrd_C(\U(\A)),
$$
where the inclusion follows from the Lemma~3 applied to the Azumaya
algebra $\A$ over the semilocal ring $C$ and extension $C_T/C$.

Hence, we can take the norm homomorphism as the transfer map and we get {\bf TE}.

Since we have the identity $\fr(T_1\times T_2)=\fr(T_1)\times \fr(T_2)$, 
additivity axiom {\bf TA} holds.

To prove {\bf TB} we consider the proof of this axiom for Linear Case 
but now the diagram ($*$) is induced by the functor 
$\R^1_{C/R}(G_m): T \mapsto \U(C_T)$ coming from
the short exact sequence of group schemes:
$$
1 \ar \R^1_{C/R}(G_m) \ar \R_{C/R}(G_m) \overset{\n^C_R}\to{\ar} G_m \ar 1,
$$
where $\R_{C/R}$ denotes Weil restriction.

We know that the functor $\R_{C/R}(G_m): T \mapsto C_T^*$
satisfies {\bf TB} and the involution is compatible with the norm map, hence, the kernel
$\R^1_{C/R}(G_m)$ satisfies axiom {\bf TB}, too.

Summarizing our discussion we have proved that the functors $\f$ and $\fr$ satisfy axioms 
{\bf C},  {\bf TE}, {\bf TA}, {\bf TB} and {\bf E}. \qed

Applying now our theorem (Non-Constant Case) to the functors $\f$ and $\fr$ 
we get the statement (a) of Unitary Case:

\proclaim{Theorem {\rm (a)}}
The map $\U(C)/\nrd(\U(\A)) \ar \U(C_K)/\nrd(\U(\A_K))$ induced 
by the canonical inclusion is injective.
\endproclaim

To show the statement (b) let look on the short exact
sequence of  smooth group schemes:
$$
1 \ar \SU_{1,\A} \ar \U_{1,\A} \overset{\nrd}\to{\ar} \R^1_{C/R}(G_m) \ar 1,
$$
where the $\U_{1,\A}: T \mapsto \U(A_T)$ for any $R$-algebra $T$.

It induces the long exact cohomology sequence:
$$
1 \ar \SU(\A) \ar \U(\A) \overset{\nrd}\to{\ar} \U(C)
\ar \hml(R,\SU_{1,\A}) \ar \hml(R,\U_{1,\A}) \ar \cdots
$$ 

Observe that the set $\hml(R,\U_{1,\A})$ is not trivial in general.

The canonical inclusion $R \hookrightarrow K$ induces the natural map
on our long cohomology sequence, so that the diagram 
$$
\matrix
1 & \lra{} & \SU(\A) & \lra{} & \U(\A) & \lra{\nrd} & \U(C) & \lra{} &  \hml(R,\SU_{1,\A}) & \lra{} & \hml(R,\U_{1,\A}) \cr
 & & \da &  & \da & & \da & & \da && \da \cr
1 & \lra{} & \SU(\A_K) & \lra{} & \U(\A_K) & \!\!\!\!\underset{\nrd_K}\to{\ar}\!\!\!\! & \U(C_K) & \lra{} &  \hml(K,\SU_{1,\A_K}) & \lra{} & \hml(K,\U_{1,\A_K}) 
\endmatrix
$$
commutes.

Taking the cokernels from the left side we get:
$$
\CD
1 @>>> \U(C)/\nrd(\U(\A)) @>>> \hml(R,\SU_{1,\A}) @>>> \hml(R,\U_{1,\A})\\
@. @VVV @VVV @VVV \\
1 @>>>  \U(C_K)/\nrd(\U(\A_K)) @>>> \hml(K,\SU_{1,\A_K}) @>>> \hml(K,\U_{1,\A_K}) 
\endCD
$$

Since the left vertical arrow is injective by (a)
 and the right vertical arrow has the trivial kernel by the main result of \cite{PO2},
we get the triviality of the kernel of the middle vertical arrow, i.e we get:

\proclaim{Theorem {\rm (b)}}
The kernel of the map $\hml(R,\SU_{1,\A}) \ar \hml(K,\SU_{1,\A_K})$ 
induced by the canonical inclusion is trivial.
\endproclaim

% ***************************************************************
% ***************************************************************
% ***************************************************************
\head{4. Norm Principle for the unitary group} \endhead

In this section we will prove the following theorem (see section 3 for the definitions): 

\proclaim{Theorem {\rm (Norm Principle for the unitary group)}}
Let $R$ be a semilocal ring with infinite residue fields of characteristic different from 2.
Let $(\A,\s)$ be an Azumaya algebra with unitary involution $\s$ over $R$.
Let $C$ be a center of $\A$, 
so $C$ is the etale quadratic $R$-algebra with the standard involution $\s$. 
Then the following equality holds: 
$$
\nrd_C(\U(\A))=\nrd_C(\A^*)^{1-\s},
$$
where $c^{1-\s}=c(c^\s)^{-1}$, for any $c\in C^*$.
\endproclaim

The constant case, i.e., when $R$ is a field, was proved by A. Merkurjev
in \cite{Me} (the elementary proof of this result can be found in \cite{BP}).

To simplify the proof we will assume that $R$ is a local ring with maximal ideal $\m$
and with infinite residue field $k$.
Indeed, one can easily extend all arguments below to the semilocal case.

The proof consists of three steps:

First, we are going to prove the inclusion $\U(\A) \subset \ainv$ (see subsection 4.1). 
Thus, taking reduced norms, we will get $\nrd(\U(\A)) \subset \nrd(\ainv)$.

Afterwards, we will show the inclusion 
$\nrd(V^{1-\s}) \subset \nrd(\U(\A))$ (see subsection 4.2), 
where $V$ is an open subset in $\A^*$, 
in some specially chosen topology.

For this we will introduce some kind of Zariski topology 
on the free $R$-module $\A$.
Then by using tricky arguments 
with the intersection of two free submodules of $\A$,
we get that almost all invertible elements have their
reduced norms in the $\nrd(\U(\A))$. 
The precise description of this procedure
is a bit technical and is contained in subsection 4.3.

We will finish with `approximation' lemma which gives us the way
to prove the inclusion $\nrd(\ainv) \subset \nrd(\U(\A))$.

\subhead\nofrills{4.1. The Proof of the Inclusion $\U(\A) \subset \ainv$}\endsubhead

Let make some remarks and definitions on the structure 
of considered rings and unitary groups.

Consider the quotients $\bar{C}=C/\m C$ and $\bar{R}=k$.
By definition $\bar{C}$ is the etale quadratic algebra over $k$
with the standard involution $\s$.
Therefore, $\bar{C}$ is the quadratic field extension $k(\sqrt{s})$, $s\in k^*$, of $k$ 
or it is the product of two fields $k\times k$.

Let $\U(\bar{C})$ be the unitary group of $\bar{C}$. 
By $\overline{\U(C)}$ we will denote the image of the unitary group $\U(C)$
under the canonical projection $C \ra C/\m C$.
Recall that 
$$
\U(C)=\{c\in C \mid \n(c)=cc^\s=1\}={C^*}^{1-\s},
$$
where $\n: C^* \ra R^*$ is the norm map and 
the last equality holds because of Hilbert~90 for local rings \cite{Se}, thus,
we have $\overline{\U(C)}=\overline{{C^*}^{1-\s}}={\bar{C}^*}{^{1-\s}}=\U(\bar{C})$.

Now we are going to prove the inclusion $\U(\A) \subset \ainv$:
The main idea here is the same as in the proof of Hilbert 90.

Let $a \in \U(\A)$. Write $b_c=c + c^\s a$, for every $c \in C^*$.
If one of these $b_c$ is invertible, then $a b_c^\s= a (c^\s + c a^\s)= b_c$, 
hence, $a= b_c {(b_c^\s)}^{-1} \in \ainv$.

Denote $\l= c{(c^\s)}^{-1} \in \U(C)$, then $b_c=c^\s(\l +a)$.
Thus, our aim is to show that for every $a \in \U(\A)$ there exists $\l \in \U(C)$ 
such that $(\l+a)\in \A^*$. 

The condition $(\l +a) \in \A^*$ is equivalent to 
$\nrd(\l + a)=\chr_a(-\l) \in C^*$,
where $\chr_a(t) \in C[t]$ is the reduced characteristic polynomial of $a$.

Since $\m C$ is the radical ideal of $C$, the element $\chr_a(-\l)$ is invertible in $C$
if and only if its quotient $\overline{\chr_a(-\l)}$ is invertible in $\bar{C}$.

Thus, our aim is to prove:
\proclaim{Lemma 4}
For every $a\in \U(\A)$ there exists $\bar{\l}\in \U(\bar{C})$ such that 
$\overline{\chr_a}(-\bar{l}) \in \bar{C}^*$, where $\overline{\chr_a}(t) \in \bar{C}[t]$
denotes the quotient of the reduced characteristic polinomial $\chr_a(t)$.
\endproclaim

\demo{Proof}
Assume that $\overline{\chr_a}(-\bar{\l})$ is non-invertible for every $\bar{\l} \in \U(\bar{C})$.

In case when $\bar{C}$ is the field 
it means that $\overline{\chr_a}(-\bar{\l})=0$ for every $\bar{\l} \in \U(\bar{C})$. 
Hence, every element $\bar{l} \in \U(\bar{C})$ is the root of the polinomial $\overline{\chr_a}(t)$.
Since the number of roots is finite we get that the unitary group $\U(\bar{C})$ must be finite.

Consider now the case when $\bar{C}=k\times k$ is the product of two fields.
Thus, the element $\bar{\l} \in \U(\bar{C})$ splits on two elements $\l'$ and $\l''$ in $k$ and
the reduced characteristic polynomial $\overline{\chr_a}(t)$ 
splits on two polynomials $\chr_a'(t)$ and $\chr_a''(t)$ over $k$ such that
$\chr_a'(\l')=0$ or $\chr_a''(\l'')=0$.

As before we conclude that either $\l'$ or $\l''$ is the root of the corresponding
polinomial and the numbers of such roots are finite.
So the unitary group $\U(\bar{C})$ consists of the points
$(x_i,y_\alpha)$ and $(x_\beta,y_i)$, where $i$ runs some finite set.

Since the involution $\s$ in this case is just the permutation map $(x,y)\mapsto (y,x)$
we get that the unitary group $\U(\bar{C})=\{\bar{c}\in \bar{C} \mid \bar{c}\bar{c}^\s=1\}$ 
consists of the points of the type $(x,x^{-1})$, $x\in k^*$.

Joining together these arguments we get that the group $\U(\bar{C})$ must be finite.

On the contrary the element of the group $\U(\bar{C})$ is just
a point on the rational quadratic curve 
$x^2- sy^2=1$ $\Longleftrightarrow$ $\n(\bar{c})=\bar{1}$
over the residue field $k$ and the number of such points is infinite. 
Thus, we get contradiction. \qed
\enddemo

\subhead\nofrills{4.2. The Proof of the Equality $\nrd({\A^*}^{1-\s})=\nrd(\U(\A))$}\endsubhead

\subsubhead{Topology on $R^n$}\endsubsubhead
We introduce topology on $R^n$ by lifting the Zariski topology given on the affine space $k^n$
via the map $R^n \xra{\can} k^n$ induced by the canonical projection.
Thus, the subbase of this topology consists of the preimages $\can^{-1}(V_f)$ of 
the main open subsets $V_f=\{ x\in k^n \mid f(x) \neq 0 \}$, 
where $f \in k[t_1,t_2,\ldots,t_n]$. 

Two important and obvious properties of this topology are:  
\  i)~every finite system of open subsets has nonempty intersection;\  and \  
ii)~any polynomial map $g:R^l \ra R^m$ is continuous. 

For instance, to see ii) it is enough to look at the commutative diagram:
$$
\CD
R^l @>g>> R^m \\
@V{\can}VV @VV{\can}V \\
k^l @>>\bar{g}> k^m
\endCD
$$
where $\bar{g}$ denotes the quotient of the polinomial map $g$. Since
the vertical arrows and the down arrow are continuous, the upper arrow
is continuous as well.

In the same way we get another important examples of continuous maps:

\itemitem{1.}~Let $A$ be an $R$-algebra and a free $R$-module of rank $n$. 
Then the regular representation 
$R^n=A \xra{rpr} \End_R A=M_n(R)=R^{n^2}$ is continuous.

\itemitem{2.}~Multiplication map 
$M_n(R) \times R^n \xra{m} R^n$, $m: (M,x) \mapsto Mx$, is continuous.

\itemitem{3.}~Let $\GL_n(R)$ be the group of invertible $R$-matrices.
Then $\GL_n(R)$ is open in $M_n(R)$ and the inverse map 
$\GL_n(R) \xra{inv} \GL_n(R)$, $inv: M \mapsto M^{-1}$,
is continuous in the induced topology on $\GL_n(R)$.

\subsubhead{The Proof of  the Inclusion of $\nrd(\ainv) \subset \nrd(\U(A))$}\endsubsubhead
Consider the Azumaya algebra $\A$ as the free R-module of rank $2m$ with
the topology constructed above.
Further we will always identify $\A$ with $R^{2m}$.

Define two free submodules of $\A$ of rank $m$:
$$
\Ap =\{ x\in \A \mid x=x^\s \} \quad \text{and} \quad \Am =\{ x\in \A \mid x=-x^\s \}.
$$
It is easy to see that $\A$ is the direct sum of the $R$-modules $\Am$ and $\Ap$.
Moreover, we have $\Am =\Ap \sqrt{b}$ and $\Ap =\Am \tfrac{\sqrt{b}}{b}$.

Consider now the intersection $a\Am \cap (R\cdot 1_{\A} \oplus \Am )$, for
some chosen $a\in \A^*$. 
We claim that (we will prove this  in the subsection 4.3) 
for almost all $a$ there exists an invertible element in
this intersection. 

The last means that there exists a non-empty open $V \subset \A^*$
such that for every $a\in V$ there exist $r\in R$, $u\in \Am$ and 
invertible element $v\in \Am$ with the property $r+u=av$.

Take $a\in V$. By the property above we may write $a=(r+u)v^{-1}$ and 
$$
\nrd(a^{1-\s})=\nrd(-(r+u)v^{-1}(r-u)^{-1}v)=
\nrd(-1)\nrd(\tfrac{r+u}{r-u}),
$$
but the elements $\tfrac{r+u}{r-u}$ and $-1$ lie in the unitary group  $\U(\A)$,
thus, we get that 
$$\nrd(a)^{1-\s} \in \nrd(\U(\A)).$$ 

The following `approximation' lemma finishes our proof:

\proclaim{Lemma 5}
Let $V \subset \A^*$ is the open subset of $\A^*$, then for every
$a\in \A^*$ there exist $v_1$, $v_2 \in V$ such that $a=v_1v_2$.
\endproclaim

\demo{Proof} By (3) the subset $V^{-1}$ is open in $\A^*$, thus, $aV^{-1}$ is open in $\A^*$. 
Since the intersection $V\cap aV^{-1}$ is nonempty, there are $v_1$, $v_2 \in V$, such that
$v_1=av_2^{-1}$ and we have proved the lemma. \qed
\enddemo

Now let $a \in \A^*$, $a=v_1v_2$,
where $v_1$, $v_2 \in V$ and
$\nrd(V)^{1-\s} \subset \nrd(\U(\A))$, then
$$
\nrd(a^{1-\s}) = \nrd(v_1)^{1-\s} \nrd(v_2)^{1-\s} \in \nrd(\U(\A)  ).
$$

This completes the proof of Equality.

\subhead\nofrills{4.3. The Proof of Existance of Invertible Element} \endsubhead

We are fixing the basis of $\A$ over $R$: 
$$\A=\Am \oplus \Ap=\{\sqrt{b}\cdot 1_{\A},\sqrt{b}e_2,\sqrt{b}e_3,\ldots,\sqrt{b}e_m\} 
\oplus 
\{1_{\A},e_2,e_3,\ldots,e_m\}.$$ 

The element $\sqrt{b}\in \Am$ in this basis is represented by the matrix
$
\left(
\smallmatrix
0_m & E_m \\
bE_m & 0_m
\endsmallmatrix
\right).
$

Consider two continuous maps $M_{2m}(R) \xra{\pr_1} M_{m-1}(R)$
and $M_{2m}(R) \xra{\pr_2} M_{m-1}(R)$, where $pr_1$ and $pr_2$ are the projections:
$$
M={(m_{ij})}_{i=1,\ldots,2m}^{j=1,\ldots,2m} \overset{\pr_1}\to\longrightarrow 
{(m_{ij})}_{i=m+2,\ldots,2m}^{j=2,\ldots,m}=N,
$$
$$
M={(m_{ij})}_{i=1,\ldots,2m}^{j=1,\ldots,2m} \overset{pr_2}\to\longrightarrow 
{(m_{ij})}_{i=m+2,\ldots,2m}^{j=1}=\bar{c}.
$$

In other words, $\pr_1$ sends any matrix $M$ to its left-down corner $(m-1)\times m$
without the first column denoted by $\bar{c}$, i.e., to the $(m-1)\times (m-1)$ matrix $N$. And $\pr_2$ sends
the matrix $M$ to the column $\bar{c}$.

Consider now continuous map 
$\psi: \A \xra{rpr} M_{2m}(R) \xra{\pr_1} M_{m-1}(R)$, where
$rpr$ is the representation of $\A$ as the $R$-module in our fixed basis.

Let $V_1$ be the preimage of the open subset 
$\GL_{m-1}(R) \subset M_{m-1}(R)$ under the map $\psi$. 

The intersection  $V_2=V_1\cap \A^*$  is the open subset in $\A^*$ 
(the subset $\A^*$ is open in $\A$).

Define the map $\omega$ as the composition:
$$
\omega: V_2 \overset{(\rho,-pr_2)}\to\longrightarrow \GL_{m-1}(R)\times R^{m-1}
\overset{m}\to\longrightarrow R^{m-1} 
\overset{i}\to\longrightarrow R^{2m}=\A,
$$
where the map $\rho$ is the composition 
$\rho: a \overset{\psi}\to\mapsto N \overset{inv}\to\mapsto N^{-1}$,
$m: (N,\bar{c}) \mapsto N\bar{c}=\bar{d}$ is the multiplication and 
$i: \bar{d}\mapsto (1,\bar{d},0,\ldots,0)^t=v$ is the inclusion.

Clearly, $\omega$ is continuous, so the preimage $V=\omega^{-1}(\A^*)$ 
of the open subset $\A^* \subset \A$ is open, too.
Since $\sqrt{b}\in V$, the subset $V$ is not empty.

We claim that for any $a\in V$
the product $av$, where $v=\omega(a)$, is the required invertible element.
Indeed, by the very construction, $v \in \Am$ and $av$ is invertible. 
Consider the product $av$ in our fixed basis:
$$
rpr(a)\omega(a)=
\left(\smallmatrix
a_{1,1}  & \ldots & a_{1,m}    & \ldots \\
\vdots   &        & \vdots     & \\
a_{m+1,1} & \ldots & a_{m+1,m} & \ldots \\
\bar{c} & & N & \ldots
\endsmallmatrix \right)
\left(\smallmatrix
1 \\
N^{-1}(-\bar{c}) \\
0 
\endsmallmatrix \right)
=
\left(\smallmatrix
u_1 \\
\vdots \\
u_{m+1} \\
\bar{c} + NN^{-1}(-\bar{c})
\endsmallmatrix \right)=
$$
$$
=\left(\smallmatrix
u_1 \\
\vdots \\
u_{m+1} \\
0
\endsmallmatrix \right) \in \Am \oplus R\cdot 1_{\A}, \quad \text{and we have finished.}
$$

% ***************************************************************
% ***************************************************************
% ***************************************************************
\head{5. The Case of a Local Regular Ring containing a Field} \endhead

Now we generalize all theorems of section~3 to the case
when $R$ is a local regular ring containing a field.
To do this we will consider in detailes only Linear Case.
Other cases can be deduced easily by following the discussion below.

Let $R$ be a local regular ring containing an infinite field of characteristic $\neq$ 2. 
Let $k$ be its residue field. Let $K$ be a field of fractions of $R$.
Let $\A$ be an Azumaya algebra over $R$.

Let $\f$ be the functor (subsection 3.2) defined on the category of $R$-algebras as:
$$
\f: T \mapsto T^*/\nrd(\A_T^*).
$$

Our aim is to show the generalized version of Theorem (a) of Linear Case (section 3):
\proclaim{Theorem {\rm (a)}}
The map $\f(R) \ra \f(K)$ induced by the canonical inclusion  is injective. 
\endproclaim

\demo{Proof} (compare with Proof of Theorem B, section~8, \cite{PO1})

By Popescu's theorem (\cite{PO1}, section~7) $R$ is a filtered direct limit of essentially smooth
local $k$-algebras, i.e, $R=\varinjlim R_i$, where $R_i$ is a local regular ring
of an affine smooth variety over $k$. 

Since an Azumaya $R$-algebra $\A$ is given by the finite number of 
generators and relations we can find the index $j$ 
such that the algebra $\A$ come from some Azumaya
algebra $\A_j$ over $R_j$, i.e. $\A=\A_j \otimes_{R_j} R$ (compare with 
the property {\bf A1} of subsection 3.1).

Fix the index $j$. 
Thus, we may assume that the functor $\f$ is given on the
category of $R_j$-algebras.
We may replace the filtered  direct system of the $R_i$ by the subsystem of all 
$R_i$ with $i \geq j$. 
Clearly we have $R=\varinjlim_{i \geq j} R_i$.
 
Let $\a \in \f(R)$ be such that its image $\a_K$ in $\f(K)$ is trivial.
Since $\f$ is continuous (subsection 3.2 axiom {\bf C}) we can find 
(as in the beginning of subsection 1.3) the localization $R_f$ of $R$ such that 
the image $\a_f$ in $\f(R_f)$ of the element $\a$ is trivial as well.

For a suitable index $k\geq j$ choose lift $f_k$ of $f$ in $R_k$.
Replacing the filtered  direct system of the $R_i$, $i\geq j$, by the subsystem of all 
$R_i$ with $i \geq k$ we still have  $R=\varinjlim_{i \geq k} R_i$. 
We put , for every $i\geq k$,
$f_i=\phi_{ik}(f_k)$ where the $\phi_{ik}: R_k \ra R_i$ are the transition homomorphisms.
It is easy to see that $\varinjlim_{i\geq k}(R_i)_{f_i} = R_f$.

By the very definition the functor $\f$ commutes
with filtered direct limits (see the proof of axiom {\bf C} in subsection 3.2), i.e.
the canonical map $\varinjlim_{i \geq k} \f(R_i) \ra \f(R)$ is an isomorphism.
Thus, we have
$$
\varinjlim_{i\geq k} \ker[ \f(R_i)\ra \f((R_i)_{f_i}) ] = \ker[\f(R) \ra \f(R_f)].
$$

By Theorem (a) of Linear Case (section 3) 
for any $i$ the map $\f(R_i) \ra \f(K_i)$ induced by the canonical inclusion
of local ring $R_i$ to its field of fractions $K_i$ is injective.
Since the inclusion $R_i \hra{} K_i$, $i\geq k$, is factorized through $(R_i)_{f_i}$ the map
$\f(R_i)\ra \f((R_i)_{f_i})$ is injective as well.

Thus, we get that the left side of the relation above is trivial. 
By the very assumption the element $\a \in \f(R)$ lies in the right side, therefore, it is trivial.
And we have finished. \qed
\enddemo

To see the generalized version of Theorem (b) of Linear Case 
one have to follow the corresponding discussion in the proof of Linear Case (subsection 3.2). 
Since all arguments there don't depend on the fact that $R$ is essentially smooth over $k$,
nothing will be changed.

\end{document}